\documentclass[12pt]{article}
\usepackage{amssymb}
\usepackage{latexsym}
\usepackage{amsthm}

\newtheorem{definition}{Definition}[section]
\newtheorem{theorem}[definition]{Theorem}
\newtheorem{lemma}[definition]{Lemma}
\newtheorem{corollary}[definition]{Corollary}
\newtheorem{remark}[definition]{Remark}

\newtheorem{note}[definition]{Note}
\newtheorem{assumption}[definition]{Assumption}
\newtheorem{proposition}[definition]{Proposition}

\newcommand{\beast}{\begin{eqnarray*}}
\newcommand{\eeast}{\end{eqnarray*}}

\def\K{\mathbb K}

\def\K{\mathbb K}
\def\fld{\mathbb K}

\begin{document}

\title{\bf Hessenberg Pairs of Linear Transformations}
\author{Ali Godjali}


\date{November 21, 2008}


\maketitle
\begin{abstract}
\noindent 
Let $\fld$ denote a field and $V$ denote a nonzero finite-dimensional vector space over $\fld$. 
We consider an ordered pair of linear transformations $A: V
\rightarrow V$ and $A^*: V \rightarrow V$ that satisfy (i)--(iii) below.
\begin{enumerate}
\item Each of $A, A^*$ is diagonalizable on $V$.
\item There exists an ordering $\lbrace V_i \rbrace_{i=0}^d$ of the eigenspaces of $A$ such that
\begin{eqnarray*}
A^* V_i \subseteq V_0 + V_1 + \ldots + V_{i+1} \qquad \qquad (0 \leq
i \leq d),
\end{eqnarray*}
where $V_{-1} = 0$, $V_{d+1}= 0$.
\item There exists an ordering $\lbrace V^*_i \rbrace_{i=0}^{\delta}$ of the
eigenspaces of $A^*$ such that
\begin{eqnarray*}
A V^*_i \subseteq V^*_0 + V^*_1 + \ldots +V^*_{i+1} \qquad \qquad (0
\leq i \leq \delta),
\end{eqnarray*}
where $V^*_{-1} = 0$, $V^*_{\delta+1}= 0$.
\end{enumerate}

\noindent We call such a pair a {\it Hessenberg pair} on $V$. In this paper we obtain some characterizations of Hessenberg pairs. We also explain how Hessenberg pairs are related to tridiagonal pairs.

\bigskip

\end{abstract}

\section{Introduction}
In \cite[Definition 1.1]{TD00} Ito, Tanabe and Terwilliger introduced the notion of a {\it tridiagonal pair} of linear transformations. Loosely speaking, this is a pair of diagonalizable linear transformations on a nonzero finite-dimensional vector space, each of which acts on the eigenspaces of the other in a certain restricted way. In \cite[Theorem 4.6]{TD00} Ito et. al. showed that a tridiagonal pair induces a certain direct sum decomposition of the underlying vector space, called the {\it split decomposition} \cite[Definition 4.1]{TD00}. In order to clarify this result, in the present paper
we introduce a generalization of a tridiagonal pair called a {\it Hessenberg pair}. Our main results are
summarized as follows. Let $V$ denote a nonzero finite-dimensional vector space, and let $(A,A^*)$ denote a pair of diagonalizable linear transformations on $V$.
We show that if $(A,A^*)$ induces a split decomposition of $V$, then $(A,A^*)$ is a Hessenberg pair on $V$. Moreover the converse holds provided that $V$ has no proper nonzero subspaces
that are invariant under each of $A$, $A^*$.

\medskip

\noindent The rest of this section contains precise statements
of our main definitions and results. We will use the following terms. Let $\fld$ denote a field and $V$ denote a nonzero finite-dimensional vector space over $\fld$. By a {\it linear transformation on $V$}, we mean a $\fld$-linear map from $V$ to $V$. Let $A$ denote a linear transformation on $V$ and
let $W$ denote a subspace of $V$. We call $W$ an {\it eigenspace} of
$A$ whenever $W\not=0$ and there exists $\theta \in \K$ such that
\begin{eqnarray*}
W=\lbrace v \in V \;\vert \;Av = \theta v\rbrace.
\end{eqnarray*}
In this case $\theta$ is called the {\it eigenvalue} of $A$
corresponding to $W$. We say $A$ is {\it diagonalizable}
whenever $V$ is spanned by the eigenspaces of $A$.

\begin{definition}
\label{def:autp} \rm By a \emph{Hessenberg pair} on $V$, we mean an
ordered pair $(A, A^*)$ of linear transformations on $V$ that
satisfy (i)--(iii) below.   

\begin{enumerate}
\item Each of $A, A^*$ is diagonalizable on $V$.
\item There exists an ordering $\lbrace V_i \rbrace_{i=0}^d$ of the eigenspaces of $A$ such that
\begin{equation}
\label{eq:astaraction} A^* V_i \subseteq V_0 + V_1 + \ldots +
V_{i+1} \qquad \qquad (0 \leq i \leq d),
\end{equation}
where $V_{-1} = 0$, $V_{d+1}= 0$.
\item There exists an ordering $\lbrace V^*_i \rbrace_{i=0}^{\delta}$ of the
eigenspaces of $A^*$ such that
\begin{equation}
\label{eq:aaction} A V^*_i \subseteq V^*_0 + V^*_1 + \ldots
+V^*_{i+1} \qquad \qquad (0 \leq i \leq \delta),
\end{equation}
where $V^*_{-1} = 0$, $V^*_{\delta+1}= 0$.
\end{enumerate}
\end{definition}

\begin{note}              \samepage
\rm
It is a common notational convention to use $A^*$ to represent the conjugate-transpose of $A$. We are not using this convention. In a Hessenberg pair $(A,A^*)$ the linear transformations $A$ and
$A^*$ are arbitrary subject to (i)--(iii) above.
\end{note}

\begin{note}
\rm The term Hessenberg comes from matrix theory. A square matrix is called
\emph{upper Hessenberg} whenever each entry below the subdiagonal
is zero \cite[p. 28]{HESS}.
\end{note}

\noindent Referring to Definition 1.1, the orderings
$\lbrace V_i \rbrace_{i=0}^d$ and $\lbrace V^*_i \rbrace_{i=0}^{\delta}$ are not
unique in general. To facilitate our discussion of these orderings we introduce some terms.
Let ($A,A^*$) denote an ordered pair of diagonalizable linear transformations on $V$. Let $\lbrace V_i \rbrace_{i=0}^d$ (resp. $\lbrace V^*_i \rbrace_{i=0}^{\delta}$) denote any ordering of the eigenspaces of $A$ (resp. $A^*$). We say that the pair ($A,A^*$) is {\it Hessenberg with respect to $(\lbrace V_i \rbrace_{i=0}^d; \lbrace V^*_i \rbrace_{i=0}^{\delta})$} whenever these orderings satisfy (\ref{eq:astaraction}) and (\ref{eq:aaction}). 
Often it is convenient to focus on eigenvalues rather than eigenspaces. 
Let $\lbrace \theta_i \rbrace_{i=0}^d$ (resp. $\lbrace \theta^*_i \rbrace_{i=0}^{\delta}$) denote the ordering of the eigenvalues of
$A$ (resp. $A^*$) that corresponds to
$\lbrace V_i \rbrace_{i=0}^d$ (resp. $\lbrace V^*_i \rbrace_{i=0}^{\delta}$).
We say that the pair ($A,A^*$) is {\it Hessenberg with respect to $(\lbrace \theta_i \rbrace_{i=0}^d;  \lbrace \theta^*_i \rbrace_{i=0}^{\delta})$} whenever ($A,A^*$) is Hessenberg with respect to ($\lbrace V_i \rbrace_{i=0}^d; \lbrace V^*_i \rbrace_{i=0}^{\delta}$).  

\begin{definition}
\label{def:irreducible} \rm Let $(A, A^*)$ denote an ordered pair of linear
transformations on $V$. We say that the pair $(A, A^*)$ is {\it irreducible}
whenever there is no subspace $W$ of $V$ such that $AW\subseteq W$,
$A^*W\subseteq W$, $W \neq 0$, $W \neq V$.
\end{definition}

\noindent We are primarily interested in the irreducible Hessenberg
pairs. However for parts of our argument the irreducibility assumption is not needed.

\smallskip

\noindent As we will see in Proposition \ref{lem:dvsdel}, for an irreducible
Hessenberg pair the scalars $d$ and $\delta$ from Definition \ref{def:autp} are equal.

%

\medskip

\noindent We now turn to the notion of a split decomposition. We will define this notion after a few preliminary comments. By a {\it decomposition} of $V$ we mean a sequence $\lbrace U_i
\rbrace_{i=0}^d$ consisting of nonzero subspaces of $V$ such that
\begin{eqnarray}
\label{eq:vdecomp} \qquad \qquad V=U_0+U_1+\cdots + U_d \qquad
\qquad (\hbox{direct sum}).
\end{eqnarray}
For notational convenience we set $U_{-1} =  0$, $U_{d+1} = 0$. For an example of a decomposition, let
$A$ denote a diagonalizable linear transformation
on $V$. Then any ordering of the eigenspaces of $A$ is a decomposition of
$V$.

\begin{lemma}
\label{lem:splitsetup} Let $A$ denote a linear transformation on
$V$. Let $\lbrace U_i \rbrace_{i=0}^d$ denote a decomposition of
$V$ and let $\lbrace \theta_i \rbrace_{i=0}^d$ denote a sequence
of mutually distinct elements of $\fld$.  Assume 
\begin{eqnarray}
\label{eq:split} (A-\theta_{i}I)U_i \subseteq U_{i+1} \qquad \qquad (0 \leq i \leq d).
\end{eqnarray}

\noindent Then $A$ is diagonalizable and $\lbrace \theta_i
\rbrace_{i=0}^d$ are the eigenvalues of $A$.
\end{lemma}

\noindent {\it Proof:} From
(\ref{eq:split}) we see that, with respect to an appropriate basis
for $V$, $A$ is represented by a lower triangular matrix which has
diagonal entries $\lbrace \theta_i \rbrace_{i=0}^d$, with $\theta_i$ appearing dim($U_i$) times for $0 \leq i \leq d$. Therefore $\lbrace
\theta_i \rbrace_{i=0}^d$ are the roots of the characteristic polynomial of $A$.
It remains to show that $A$ is diagonalizable. From  (\ref{eq:split})
we see that $\prod_{i=0}^{d} (A- \theta_iI)$ vanishes on $V$. By
this and since $\lbrace \theta_i \rbrace_{i=0}^d$ are distinct we
see that the minimal polynomial of $A$ has distinct roots. Therefore $A$ is diagonalizable and the result follows. \hfill $\Box$

\medskip


\begin{definition}
\label{def:sp}
\rm

Let $d$ denote a nonnegative integer. Let $A$ (resp. $A^*$) denote a
diagonalizable linear transformation on $V$ with eigenvalues
$\lbrace \theta_i \rbrace_{i=0}^d$ (resp. $\lbrace \theta^*_i
\rbrace_{i=0}^d$). By an {\it $(A, A^*)$-split decomposition of $V$ with respect to $(\lbrace \theta_i\rbrace_{i=0}^d; \lbrace
\theta^*_i\rbrace_{i=0}^{d})$}, we mean a decomposition $\lbrace U_i
\rbrace_{i=0}^d$ of $V$ such that both 
\begin{eqnarray}
&(A-\theta_{d-i}I)U_i \subseteq U_{i+1},& \label{eq:split1}
\\
&(A^*-\theta^*_iI)U_i \subseteq U_{i-1}& \label{eq:split2}
\end{eqnarray}
for $0 \leq i \leq d$. 
\end{definition}

\noindent As we will see in Corollary \ref{cor:uniqueness}, the $(A, A^*)$-split decomposition of $V$ with respect to \\
$( \lbrace \theta_i\rbrace_{i=0}^d; \lbrace
\theta^*_i\rbrace_{i=0}^{d})$ is unique if it exists.

\bigskip

\noindent The main results of this paper are the following two theorems and subsequent corollary. 

\begin{theorem}
\label{thm:hesstosplit} Let $d$ denote a nonnegative integer. Let
$A$ (resp. $A^*$) denote a diagonalizable linear transformation on
$V$ with eigenvalues $\lbrace \theta_i \rbrace_{i=0}^d$ (resp. $\lbrace \theta^*_i \rbrace_{i=0}^d$). Suppose that the pair ($A,A^*$) is irreducible, and Hessenberg with respect to $(\lbrace \theta_i \rbrace_{i=0}^d; \lbrace \theta^*_i \rbrace_{i=0}^{d})$. Then there exists
an $(A, A^*)$-split decomposition of $V$ with respect to $(\lbrace
\theta_{i} \rbrace_{i=0}^d; \lbrace \theta^*_i \rbrace_{i=0}^{d})$.
\end{theorem}

\begin{theorem}
\label{thm:splittohess} Let $d$ denote a nonnegative integer. Let
$A$ (resp. $A^*$) denote a diagonalizable linear transformation on
$V$ with eigenvalues $\lbrace \theta_i \rbrace_{i=0}^d$ (resp. $\lbrace \theta^*_i \rbrace_{i=0}^d$). Suppose that there exists an
$(A, A^*)$-split decomposition of $V$ with respect to $(\lbrace
\theta_{i} \rbrace_{i=0}^d; \lbrace \theta^*_i\rbrace_{i=0}^{d})$.
Then the pair $(A, A^*)$ is Hessenberg with respect to $(\lbrace \theta_i
\rbrace_{i=0}^d$; $\lbrace \theta^*_i \rbrace_{i=0}^d)$.
\end{theorem}

\noindent Combining Theorem \ref{thm:hesstosplit} and Theorem \ref{thm:splittohess} we obtain the following
corollary.

\begin{corollary}
\label{cor:mainresult} Let $d$ denote a nonnegative integer. Let $A$
(resp. $A^*$) denote a diagonalizable linear transformation on $V$
with eigenvalues $\lbrace \theta_i \rbrace_{i=0}^d$ (resp. $\lbrace
\theta^*_i \rbrace_{i=0}^d$). Assume the pair $(A, A^*)$ is irreducible. Then the
following (i), (ii) are equivalent.
\begin{enumerate}
\item
The pair $(A, A^*)$ is Hessenberg with respect to $(\lbrace \theta_i
\rbrace_{i=0}^d$; $\lbrace \theta^*_i \rbrace_{i=0}^d)$.
\item
There exists an $(A, A^*)$-split decomposition of $V$ with respect
to $(\lbrace \theta_{i} \rbrace_{i=0}^d; \lbrace
\theta^*_i\rbrace_{i=0}^{d})$.
\end{enumerate}
\end{corollary}

\section{The Proof of Theorem \ref{thm:hesstosplit}}

In this section we give a proof of Theorem \ref{thm:hesstosplit}. Along the way, we show that the scalars $d$
and $\delta$ from Definition \ref{def:autp} are equal. We will refer to the following setup. 

\begin{assumption}
\label{def:vij} \rm Let $A$ (resp. $A^*$) denote a diagonalizable linear transformation on
$V$ with eigenvalues $\lbrace \theta_i \rbrace_{i=0}^d$ (resp. $\lbrace \theta^*_i \rbrace_{i=0}^{\delta}$).
Let $\lbrace V_i \rbrace_{i=0}^d$ (resp. $\lbrace V^*_i \rbrace_{i=0}^{\delta}$) denote the corresponding eigenspaces of $A$ (resp. $A^*$).
We assume that the pair ($A,A^*$) is irreducible and Hessenberg with respect to ($\lbrace \theta_i \rbrace_{i=0}^d; \lbrace \theta^*_i \rbrace_{i=0}^{\delta}$). For all integers $i$ and $j$ we set
\begin{eqnarray}
V_{ij} =
(V_0+\cdots + V_i)\cap (V^*_{0}+\cdots + V^*_{j}). \label{eq:vij}
\end{eqnarray}
We interpret the sum on the left in (\ref{eq:vij}) to be 0 (resp. $V$) if $i<0$  (resp. $i>d$). We interpret the sum on the right in
(\ref{eq:vij}) to be 0 (resp. $V$) if $j<0$  (resp. $j>\delta$).
\end{assumption}

\begin{lemma}
\label{lem:basic} With reference to Assumption \ref{def:vij}, the
following (i), (ii) hold for $0 \leq i\leq d$ and $0 \leq j\leq
\delta$.
\begin{enumerate}
\item
$V_{i\delta} =
 V_0+\cdots+V_i$.
\item
$V_{dj} =
 V^*_{0}+\cdots+V^*_{j}$.
\end{enumerate}
\end{lemma}
\noindent {\it Proof:} (i) Set $j=\delta$ in (\ref{eq:vij}) and use
the fact that $V = V^*_{0} + \cdots + V^{*}_{\delta}$.
\\
(ii) Set $i=d$ in (\ref{eq:vij}) and use the fact that $V = V_0 +
\cdots + V_d$. \hfill $\Box $
\medskip

\begin{lemma}
\label{lem:kact} With reference to Assumption \ref{def:vij}, the
following (i), (ii) hold for $0 \leq i\leq d$ and $0 \leq j\leq
\delta$.
\begin{enumerate}
\item
$(A - \theta_iI)V_{ij} \subseteq V_{i-1,j+1}$.
\item
$(A^*- \theta^*_jI)V_{ij} \subseteq V_{i+1,j-1}$.
\end{enumerate}
\end{lemma}
\noindent {\it Proof:} (i) Since $V_i$ is the eigenspace of $A$
corresponding to the eigenvalue $\theta_i$, we have
\begin{eqnarray}
(A-\theta_iI)\sum_{h=0}^i V_h = \sum_{h=0}^{i-1} V_h.
\label{eq:kinv1}
\end{eqnarray}
Using (\ref{eq:aaction})  we find
\begin{eqnarray}
(A - \theta_iI)\sum_{h=0}^j V^*_h \subseteq \sum_{h=0}^{j+1} V^*_h.
\label{eq:kinv2}
\end{eqnarray}
Evaluating $(A - \theta_iI)V_{ij}$ using
(\ref{eq:vij})--(\ref{eq:kinv2}), we find it is contained in
$V_{i-1,j+1}$.
\\
\noindent (ii) Using (\ref{eq:astaraction}) we find
\begin{eqnarray}
(A^*-\theta^*_jI)\sum_{h=0}^i V_h \subseteq \sum_{h=0}^{i+1} V_h.
\label{eq:ksinv1}
\end{eqnarray}
Since $V^*_j$ is the eigenspace of $A^*$ corresponding to the eigenvalue
$\theta^{*}_{j}$, we have
\begin{eqnarray}
(A^*-\theta^*_jI)\sum_{h=0}^j V^*_h =\sum_{h=0}^{j-1} V^*_h.
\label{eq:ksinv2}
\end{eqnarray}

\noindent Evaluating $(A^*- \theta^*_jI)V_{ij}$ using
(\ref{eq:vij}), (\ref{eq:ksinv1}), (\ref{eq:ksinv2}), we find it
is contained in $V_{i+1,j-1}$. \hfill $\Box $
\medskip

\begin{proposition}
\label{lem:dvsdel} With reference to Assumption \ref{def:vij}, the
scalars $d$ and $\delta$ from Definition \ref{def:autp} are
equal. Moreover,
\begin{eqnarray}
\label{eq:ijd}
V_{ij}= 0 \quad {\rm if}\quad i+j<d,\qquad  \qquad (0
\leq i,j\leq d).
\end{eqnarray}
\end{proposition}

\noindent {\it Proof:} For all nonnegative integers $r$ such that
$r\leq d$ and $r \leq \delta$, we define
\begin{eqnarray}
W_r=V_{0r}+V_{1,r-1}+\cdots +V_{r0}. \label{eq:wr}
\end{eqnarray}
We have $AW_r \subseteq W_r$ by Lemma \ref{lem:kact}(i) and $A^*W_r
\subseteq W_r$ by Lemma \ref{lem:kact}(ii). Now $W_r=0$ or $W_r=V$
since the pair $(A, A^*)$ is irreducible. Suppose for the moment that
$r\leq d-1$. Each term on the right in (\ref{eq:wr}) is contained in
$V_0+\cdots + V_r$ so $W_r \subseteq V_0+\cdots + V_r$. Thus
$W_r\not=V$ and hence $W_r=0$. Next suppose $r=d$. Then
$V_{d0}\subseteq W_r$. Recall $V_{d0}=V^*_0$ by Lemma
\ref{lem:basic}(ii) and $V^*_0\not=0$ so $V_{d0}\not=0$. Now
$W_r\not=0$ so $W_r=V$. We have now shown that $W_r=0$ if $r\leq
d-1$ and $W_r=V$ if $r=d$. Similarly $W_r=0$ if $r\leq \delta-1$ and
$W_r=V$ if $r=\delta$. Now $d=\delta$; otherwise we take
$r=\mbox{min}(d,\delta)$ in our above comments and find $W_r$ is
both $0$ and $V$, for a contradiction. The result follows. \hfill
$\Box $

\medskip

\begin{lemma}
\label{lem:vidissplit} With reference to Assumption \ref{def:vij}, the sequence
$\lbrace V_{d-i,i} \rbrace_{i=0}^d$ is an $(A, A^*)$-split
decomposition of $V$ with respect to $(\lbrace \theta_{i}
\rbrace_{i=0}^d; \lbrace \theta^*_i \rbrace_{i=0}^{d})$.
\end{lemma}

\noindent {\it Proof:} Observe that (\ref{eq:split1}) follows from
Lemma \ref{lem:kact}(i) and (\ref{eq:split2}) follows from Lemma
\ref{lem:kact}(ii). It remains to show that the sequence $\lbrace V_{d-i,i} \rbrace_{i=0}^d$ is a decomposition. We first
show
\begin{eqnarray}
V= \sum_{i=0}^d V_{d-i,i}. \label{eq:usumv}
\end{eqnarray}
Let $W$ denote the sum on the right in (\ref{eq:usumv}). We have
$AW\subseteq W$ by Lemma \ref{lem:kact}(i) and   $A^*W\subseteq
W$ by Lemma \ref{lem:kact}(ii). Now $W=0$ or $W=V$ by the
irreducibility assumption. Observe that $W$ contains $V_{d0}$
and $V_{d0}=V^*_0$ is nonzero so $W\not=0$. We conclude that $W=V$
and (\ref{eq:usumv}) follows. Next we show that the sum
(\ref{eq:usumv}) is direct. To do this we show that
\begin{eqnarray}
\label{expvdi}
 (V_{d0}+V_{d-1,1}+\cdots +V_{d-i+1,i-1})\cap V_{d-i,i}
\end{eqnarray}
is zero for $1 \leq i \leq d$. Let $i$ be given. From the construction
\begin{eqnarray*}
V_{d-j,j} \subseteq V^*_{0} + V^*_{1} + \cdots + V^*_{i-1}
\end{eqnarray*}
for $0 \leq j \leq i-1$, and
\begin{eqnarray*}
V_{d-i,i}\subseteq V_0 + V_1 + \cdots + V_{d-i}.
\end{eqnarray*}
Therefore (\ref{expvdi}) is contained in 
\begin{eqnarray}
\label{expr}
(V_0 + V_1 + \cdots + V_{d-i})\cap
(V^*_{0} + V^*_{1} + \cdots + V^*_{i-1}).
\end{eqnarray}
But (\ref{expr}) is equal to  $V_{d-i,i-1}$ and this is zero by (\ref{eq:ijd}), so (\ref{expvdi}) is zero. We have shown that the sum
(\ref{eq:usumv}) is direct. Next we show that $V_{d-i,i} \neq 0$ for $0 \leq i \leq d$. Suppose there exists an
integer $i$ $(0 \leq i \leq d)$ such that $V_{d-i,i}=0$. Observe
that $i\not=0$ since $V_{d0}=V^*_{0}$ is nonzero and $i\not=d$
since $V_{0d}=V_0$ is nonzero. Set
\begin{eqnarray*}
W = V_{d0}+V_{d-1,1} + \cdots+V_{d-i+1,i-1}
\end{eqnarray*}
and observe that $W\not=0$ and $W\not=V$ by our above remarks. By
Lemma \ref{lem:kact}(ii), we find $A^*W\subseteq W$. By Lemma
\ref{lem:kact}(i) and since $V_{d-i,i}=0$, we find $AW\subseteq W$.
Now $W=0$ or $W=V$ by our irreducibility assumption, which yields a
contradiction. We conclude that $V_{d-i,i} \neq 0$ for $0 \leq i \leq d$. We have now shown that the sequence $\lbrace
V_{d-i,i} \rbrace_{i=0}^d$ is a decomposition of $V$ and we are
done. \hfill $\Box$

\medskip

\noindent Theorem \ref{thm:hesstosplit} is immediate from Lemma \ref{lem:vidissplit}.

\section{The Proof of Theorem \ref{thm:splittohess}}

In this section we give a proof of Theorem \ref{thm:splittohess}. Along the way, we show that the split decomposition from Definition \ref{def:sp} is unique if it exists. The following
assumption sets the stage.

\begin{assumption}
\label{ass:splitexist}
\rm
Let $d$ denote a nonnegative integer. Let
$A$ (resp. $A^*$) denote a diagonalizable linear transformation on
$V$ with eigenvalues $\lbrace \theta_i \rbrace_{i=0}^d$ (resp. $\lbrace \theta^*_i \rbrace_{i=0}^d$). Let $\lbrace V_i
\rbrace_{i=0}^d$ (resp. $\lbrace V^*_i \rbrace_{i=0}^d$) denote the
corresponding eigenspaces of $A$ (resp. $A^*$). We
assume that there exists a decomposition $\lbrace U_i
\rbrace_{i=0}^d$ of $V$ that is $(A, A^*)$-split with respect to $(
\lbrace \theta_i\rbrace_{i=0}^d; \lbrace \theta^*_i\rbrace_{i=0}^{d}
)$.
\end{assumption}

\begin{lemma}
\label{thm:maincharls}
With reference to Assumption \ref{ass:splitexist}, for $0 \leq i \leq d$ both
\begin{eqnarray}
U_{i}+U_{i+1}+\cdots +U_{d} &=& V_{0}+V_{1} +\cdots + V_{d-i},
\label{eq:vsumitod}
\\
U_0+U_1+\cdots +U_i &=& V^*_{0}+V^*_{1} +\cdots + V^*_{i}.
\label{eq:vsumzeroi}
\end{eqnarray}
\end{lemma}

\noindent {\it Proof:} First consider (\ref{eq:vsumitod}). We abbreviate
\beast
W=U_{i}+U_{i+1}+\cdots + U_d, \qquad \quad
Z=V_0+V_{1}+\cdots + V_{d-i}.
\eeast
We show $W=Z$. To obtain $Z\subseteq W$, set $X=\prod_{h=0}^{i-1} (A-\theta_{d-h}I)$, and observe
$Z=XV$ by elementary linear algebra. Using
(\ref{eq:split1}), we find $XU_j\subseteq W$ for $0 \leq j \leq d$, so $XV\subseteq W$ in view of
(\ref{eq:vdecomp}).
We now have $Z\subseteq W$. To obtain $W\subseteq Z$, set $Y=\prod_{h=i}^{d} (A-\theta_{d-h}I)$, and observe
\begin{eqnarray}
Z&=&\lbrace v \in V \;|\;Yv = 0\rbrace .
\label{eq:Keraction}
\end{eqnarray}
Using (\ref{eq:split1}), we find $YU_j=0$ for $ i \leq j \leq d$, so $YW=0$. Combining
this with (\ref{eq:Keraction}), we find $W\subseteq Z$. We now have $Z=W$ and hence
(\ref{eq:vsumitod}) holds.
Line (\ref{eq:vsumzeroi})  is similarly obtained using (\ref{eq:split2}).
\hfill $\Box $

\begin{lemma}
\label{lem:whatisui}
With reference to Assumption \ref{ass:splitexist},
\begin{eqnarray}
\label{eq:uniqueness}
U_i = (V^*_{0}+V^*_{1}+\cdots + V^*_{i})\cap
(V_{0}+V_{1}+\cdots + V_{d-i}) \qquad (0 \leq i \leq d).
\end{eqnarray}
\end{lemma}

\noindent {\it Proof:}  Since $\lbrace U_i \rbrace_{i=0}^d$ is a
decomposition of $V$,
\begin{eqnarray}
\label{uni}
U_i &=& (U_0+U_1+\cdots + U_i)\cap (U_i+U_{i+1}+\cdots +
U_d) \qquad (0 \leq i \leq d).
\end{eqnarray}
Evaluating (\ref{uni}) using (\ref{eq:vsumitod}),
(\ref{eq:vsumzeroi}) we obtain (\ref{eq:uniqueness}).
\hfill $\Box$

\begin{corollary}
\label{cor:uniqueness}
With reference to Assumption \ref{ass:splitexist}, the split decomposition $\lbrace U_i \rbrace_{i=0}^d$ is uniquely determined by the given orderings of the eigenvalues $\lbrace \theta_i \rbrace_{i=0}^d$ and
$\lbrace \theta^*_i \rbrace_{i=0}^d$. 
\end{corollary}

\noindent {\it Proof:} Immediate from Lemma \ref{lem:whatisui}. \hfill $\Box$

\medskip

\begin{lemma}
\label{lem:hessmotive}
With reference to Assumption \ref{ass:splitexist}, for $0 \leq i \leq d$ both
\begin{eqnarray}
\label{inverseorder}
A^* V_i \subseteq V_{0} + V_1 + \cdots + V_{i+1},
\\
A V^*_i \subseteq V^*_0 + V^*_1 + \cdots + V^*_{i+1}.
\label{sameorder}
\end{eqnarray}
Moreover $(A, A^*)$ is a Hessenberg pair on $V$.
\end{lemma}

\noindent {\it Proof:} To obtain (\ref{inverseorder}), observe
\begin{eqnarray*}
A^*V_{i} &\subseteq& A^* \sum_{h=0}^{i} V_{h}
\\
&=& A^* \sum_{h=d-i}^{d} U_{h} \qquad \qquad (\hbox{by (\ref{eq:vsumitod})})
\\
&\subseteq& \sum_{h=d-i-1}^{d} U_{h} \qquad \qquad \ \ (\hbox{by (\ref{eq:split2})})
\\
&=& \sum_{h=0}^{i+1} V_{h} \qquad \qquad \qquad  \  (\hbox{by (\ref{eq:vsumitod})}).
\end{eqnarray*}

To obtain (\ref{sameorder}), observe
\begin{eqnarray*}
AV^*_{i} &\subseteq& A \sum_{h=0}^{i} V^*_{h}
\\
&=& A \sum_{h=0}^{i} U_{h} \qquad \qquad \ \ \ \ \  (\hbox{by (\ref{eq:vsumzeroi})})
\\
&\subseteq& \sum_{h=0}^{i+1} U_{h} \qquad \qquad \ \ \ \ \ \ \ \ (\hbox{by (\ref{eq:split1})})
\\
&=& \sum_{h=0}^{i+1} V^*_{h} \qquad \qquad \ \ \ \ \ \ \ \ (\hbox{by (\ref{eq:vsumzeroi})}).
\end{eqnarray*}
\hfill $\Box$

\medskip

\noindent Theorem \ref{thm:splittohess} is immediate from Lemma \ref{lem:hessmotive}. 

\bigskip

\noindent We finish this section with a comment.
\begin{corollary}
With reference to Assumption \ref{ass:splitexist}, for $0 \leq i \leq d$ the dimensions of $V_{d-i}$, $V^*_i$, $U_i$ are the same.
\end{corollary}

\noindent {\it Proof:} Recall that $\lbrace V_i \rbrace_{i=0}^d$ and
$\lbrace U_i \rbrace_{i=0}^d$ are decompositions of $V$. By this and 
(\ref{eq:vsumitod}), 
\begin{eqnarray*}
\hbox{\rm dim}(U_{i})+ {\rm dim} (U_{i+1})+ \cdots + {\rm dim}(U_{d}) = {\rm dim}(V_{0})+ {\rm dim}(V_{1})+\cdots + {\rm dim}(V_{d-i})
\end{eqnarray*}
for $0 \leq i \leq d$. Consequently, the dimensions of $V_{d-i}$ and
$U_i$ are the same for $0 \leq i \leq d$. A similar argument using
(\ref{eq:vsumzeroi}) shows that the dimensions of $V^*_i$ and $U_i$
are the same for $0 \leq i \leq d$. The result follows. \hfill
$\Box$

\section{Hessenberg pairs and tridiagonal pairs}

In this section, we explain how Hessenberg pairs are related to tridiagonal pairs. Using this relationship we show that some results \cite[Lemma 4.5]{TD00}, \cite[Theorem 4.6]{TD00} about tridiagonal pairs are direct consequences of our results on Hessenberg pairs. We start by recalling the definition of a tridiagonal pair. 

\begin{definition}
{\rm \cite[Definition 1.1]{TD00}}
\label{def:tdpair}
\rm
By a \emph{tridiagonal pair} on $V$, we mean an ordered pair $(A, A^*)$ of linear transformations on $V$ that satisfy (i)--(iv) below.

\begin{enumerate}
\item Each of $A, A^*$ is diagonalizable on $V$.
\item There exists an ordering $\lbrace V_i
\rbrace_{i=0}^d$ of the
eigenspaces of $A$ such that
\begin{equation}
A^* V_i \subseteq V_{i-1} + V_i+ V_{i+1} \qquad \qquad (0 \leq i \leq d),
\label{eq:astaractionthreeint}
\end{equation}
where $V_{-1} = 0$, $V_{d+1}= 0$.
\item There exists an ordering $\lbrace V^*_i
\rbrace_{i=0}^{\delta}$ of
the
eigenspaces of $A^*$ such that
\begin{equation}
A V^*_i \subseteq V^*_{i-1} + V^*_i+ V^*_{i+1} \qquad \qquad (0 \leq i \leq \delta),
\label{eq:aactionthreeint}
\end{equation}
where $V^*_{-1} = 0$, $V^*_{\delta+1}= 0$.
\item The pair $(A, A^*)$ is irreducible in the sense of Definition \ref{def:irreducible}.  
\end{enumerate}
\end{definition}

\begin{definition}
\label{def:tdpairsorderings}
\rm
Let $(A, A^*)$ denote an ordered pair of diagonalizable linear transformations on $V$. Let $\lbrace V_i \rbrace_{i=0}^d$ (resp. $\lbrace V^*_i \rbrace_{i=0}^{\delta}$) denote any ordering of the eigenspaces of $A$ (resp. $A^*$). We say that the pair ($A,A^*$) is {\it tridiagonal with respect to $(\lbrace V_i \rbrace_{i=0}^d; \lbrace V^*_i \rbrace_{i=0}^{\delta})$} whenever conditions (ii)--(iv) in Definition \ref{def:tdpair} are satisfied. 
\end{definition}

\begin{remark}
\rm
With reference to Definition \ref{def:tdpairsorderings}, assume that ($A,A^*$) is tridiagonal with respect to ($\lbrace V_i \rbrace_{i=0}^d; \lbrace V^*_i \rbrace_{i=0}^{\delta}$). By \cite[Lemma 2.4]{TD00} the pair ($A,A^*$) is tridiagonal with respect to each of ($\lbrace V_{d-i} \rbrace_{i=0}^d; \lbrace V^*_i \rbrace_{i=0}^{\delta}$), ($\lbrace V_{i} \rbrace_{i=0}^d; \lbrace V^*_{\delta-i} \rbrace_{i=0}^{\delta}$), ($\lbrace V_{d-i} \rbrace_{i=0}^d; \lbrace V^*_{\delta-i} \rbrace_{i=0}^{\delta}$) and no further orderings of the eigenspaces.
\end{remark}

\noindent Hessenberg pairs and tridiagonal pairs are related as follows. 

\begin{proposition}
\label{lem:tdandaut}
Let $A$ (resp. $A^*$) denote a diagonalizable linear transformation on
$V$ with eigenspaces $\lbrace V_i \rbrace_{i=0}^d$ (resp. $\lbrace V^*_i \rbrace_{i=0}^{\delta}$). Then the following (i)--(iv) are equivalent.

\begin{enumerate}
\item
The pair $(A, A^*)$ is tridiagonal with respect to $(\lbrace V_i \rbrace_{i=0}^d; \lbrace V^*_i \rbrace_{i=0}^{\delta})$. 
\item
The pair $(A, A^*)$ is irreducible, and Hessenberg with respect to each of $(\lbrace V_i \rbrace_{i=0}^d; \lbrace V^*_i \rbrace_{i=0}^{\delta})$, $(\lbrace V_{d-i} \rbrace_{i=0}^d; \lbrace V^*_i \rbrace_{i=0}^{\delta})$, $(\lbrace V_i \rbrace_{i=0}^d; \lbrace V^*_{\delta-i} \rbrace_{i=0}^{\delta})$,   $(\lbrace V_{d-i} \rbrace_{i=0}^d; \lbrace V^*_{\delta-i} \rbrace_{i=0}^{\delta})$.
\item
The pair $(A, A^*)$ is irreducible, and Hessenberg with respect to each of $(\lbrace V_i \rbrace_{i=0}^d; \lbrace V^*_i \rbrace_{i=0}^{\delta})$, $(\lbrace V_{d-i} \rbrace_{i=0}^d; \lbrace V^*_{\delta-i} \rbrace_{i=0}^{\delta})$.
\item
The pair $(A, A^*)$ is irreducible, and Hessenberg with respect to each of $(\lbrace V_{d-i} \rbrace_{i=0}^d; \lbrace V^*_i \rbrace_{i=0}^{\delta})$, $(\lbrace V_i \rbrace_{i=0}^d; \lbrace V^*_{\delta-i} \rbrace_{i=0}^{\delta})$.
\end{enumerate}
\end{proposition}

\noindent {\it Proof:} Observe that $\lbrace V_i \rbrace_{i=0}^d$ satisfies (\ref{eq:astaractionthreeint}) if and only if both $\lbrace V_i \rbrace_{i=0}^d$ and $\lbrace V_{d-i} \rbrace_{i=0}^d$ satisfy (\ref{eq:astaraction}).
Similarly $\lbrace V^*_i \rbrace_{i=0}^{\delta}$ satisfies (\ref{eq:aactionthreeint}) if and only if
both $\lbrace V^*_i \rbrace_{i=0}^{\delta}$ and $\lbrace V^*_{d-i} \rbrace_{i=0}^{\delta}$ satisfy (\ref{eq:aaction}). The result follows. 
\hfill $\Box$

\medskip
\noindent In Proposition \ref{lem:tdandaut} we showed how Hessenberg pairs are related to tridiagonal pairs. We now use this relationship to obtain some results on tridiagonal pairs. 

\begin{theorem}
{\rm \cite[Lemma 4.5]{TD00}}
Let $(A, A^*)$ denote a tridiagonal pair as in Definition \ref{def:tdpair}. Then the scalars $d$ 
and $\delta$ from that definition are equal.
\end{theorem}

\noindent {\it Proof:} Combine Proposition \ref{lem:dvsdel} and Proposition \ref{lem:tdandaut}. \hfill $\Box$

\medskip

\begin{definition}
\rm
Let $(A, A^*)$ denote an ordered pair of diagonalizable linear transformations on $V$. Let $\lbrace \theta_i \rbrace_{i=0}^d$ (resp. $\lbrace \theta^*_i \rbrace_{i=0}^{\delta}$) denote any ordering of the eigenvalues of $A$ (resp. $A^*$). Let $\lbrace V_i \rbrace_{i=0}^d$ (resp. $\lbrace V^*_i \rbrace_{i=0}^{\delta}$) denote the corresponding ordering of the eigenspaces of $A$ (resp. $A^*$).
We say that the pair ($A,A^*$) is {\it tridiagonal with respect to $(\lbrace \theta_i \rbrace_{i=0}^d;  \lbrace \theta^*_i \rbrace_{i=0}^{\delta})$} whenever ($A,A^*$) is tridiagonal with respect to $(\lbrace V_i \rbrace_{i=0}^d; \lbrace V^*_i \rbrace_{i=0}^{\delta})$.  
\end{definition}

\begin{theorem}
{\rm \cite[Theorem 4.6]{TD00}}
\label{thm:characterization}
Let $d$ denote a nonnegative integer. Let
$A$ (resp. $A^*$) denote a diagonalizable linear transformation on
$V$ with eigenvalues $\lbrace \theta_i \rbrace_{i=0}^d$ (resp. $\lbrace \theta^*_i \rbrace_{i=0}^d$). Then the following (i)--(iv) are equivalent.

\begin{enumerate}
\item
The pair $(A, A^*)$ is tridiagonal with respect to $(\lbrace \theta_i \rbrace_{i=0}^d; \lbrace \theta^*_i \rbrace_{i=0}^d)$.
\item
The pair $(A, A^*)$ is irreducible, and there exist $(A, A^*)$-split decompositions of $V$ with respect to each of $(\lbrace \theta_{i} \rbrace_{i=0}^d; \lbrace \theta^*_i\rbrace_{i=0}^{d})$, $(\lbrace \theta_{d-i} \rbrace_{i=0}^d; \lbrace \theta^*_{i} \rbrace_{i=0}^{d})$, $(\lbrace \theta_i \rbrace_{i=0}^d; \lbrace \theta^*_{d-i} \rbrace_{i=0}^{d})$, \\ $(\lbrace \theta_{d-i} \rbrace_{i=0}^{d}; \lbrace \theta^*_{d-i} \rbrace_{i=0}^{d})$.
\item
The pair $(A, A^*)$ is irreducible, and there exist $(A, A^*)$-split decompositions of $V$ with respect to each of $(\lbrace \theta_{i} \rbrace_{i=0}^d; \lbrace \theta^*_i \rbrace_{i=0}^{d})$, $(\lbrace \theta_{d-i} \rbrace_{i=0}^d; \lbrace \theta^*_{d-i} \rbrace_{i=0}^{d})$.
\item
The pair $(A, A^*)$ is irreducible, and there exist $(A, A^*)$-split decompositions of $V$ with respect to each of $(\lbrace \theta_{d-i} \rbrace_{i=0}^d; \lbrace \theta^*_i \rbrace_{i=0}^{d})$, $(\lbrace \theta_{i} \rbrace_{i=0}^d; \lbrace \theta^*_{d-i} \rbrace_{i=0}^{d})$.
\end{enumerate}
\end{theorem}

\noindent {\it Proof:} Combine Corollary \ref{cor:mainresult} and Proposition \ref{lem:tdandaut}. \hfill $\Box$

\medskip

\section{Acknowledgement}
This paper was written while the author was a graduate student at the University of Wisconsin-Madison. The author would like to thank his advisor Paul Terwilliger for his many valuable ideas and suggestions.

\noindent Ali Godjali \hfil\break
\noindent Department of Mathematics \hfil\break
\noindent University of Wisconsin \hfil\break
\noindent Van Vleck Hall \hfil\break
\noindent 480 Lincoln Drive \hfil\break
\noindent Madison, WI 53706-1388 USA \hfil\break
\noindent email: {\tt godjali@math.wisc.edu }\hfil\break


\begin{thebibliography}{10}

\bibitem{TD00}
T.~Ito, K.~Tanabe, and P.~Terwilliger.
\newblock Some algebra related to ${P}$- and ${Q}$-polynomial association
  schemes,  in:
\newblock {\em Codes and Association Schemes (Piscataway NJ, 1999)}, Amer.
Math. Soc., Providence RI, 2001, pp.
     167--192;
{\tt arXiv:math.CO/0406556}.

\bibitem{HESS}
R. ~Horn and C. ~Johnson,
\newblock {\em Matrix Analysis}, 
\newblock Cambridge University Press, New York, 1985.
\newblock


%
%
%
%
%
%
%
%
%
%
%
%
%
%
%
%
%


















\end{thebibliography}
\end{document}